%% file: fixed-shore-journal.tex
\documentclass[12pt]{article}
\newcommand{\runningtitle}{Quantitative group testing with a fixed shore}
\input{layout}
\input{width-changes}
\hypersetup{pdftitle={\papername},pdfauthor={Fedor Karpelevitch}}
\begin{document}
\thispagestyle{plain}
\begin{center}
{\LARGE\papername\par}\vspace{.6em}
\papersignature
\end{center}
\begin{abstract}
Two disjoint sets each contain one defective item. A test reports how many
defectives lie in a chosen subset, and its choice may depend on earlier
answers. We determine the largest searchable second set for every test budget
when the first set has six, seven or eight items, extending the known answers
through five. For six items, every optimal first test splits them four and
two once at least nine tests are available. For eight, the exact formula
gains one item over a quadratic correction every three tests. We explain
the rounding in an optimal strategy and prove a matching bound over all
strategies using exact rational certificates. For an arbitrary fixed first
set, we determine the degree and leading coefficient of the loss from the
binary bound.
\end{abstract}

\section{How should the first test split six items?}
Two disjoint sets each contain one defective item. We may test any subset
of their union and learn how many defectives it contains. If the first set
has \(m\) items, how large can the other set be if we must find both defectives
in at most \(k\) tests? Write \(n(k,m)\) for this largest size, taking it to be
zero if no nonempty second set qualifies. Even if the first defective were
revealed, a test would give only a yes or a no about the second. Thus
\(n(k,m)\le 2^k\).

Take six items in the first set and ten tests. If each test reports only
whether it contains a defective, Chang and Hwang's theorem \cite{CH81}
gives an exact maximum of \(\lfloor 2^{10}/6\rfloor=170\) items in the other
set. Reporting how many defectives the test contains raises that maximum
to 973. One optimal first test selects four of the six items and 496 of the 973.
Figure~\ref{fig:first} shows all three possible answers. An answer of two
leaves four items paired with 496. Aigner's four-item formula gives
\(n(9,4)=496\), so even one more item in that block would be too many.
An answer of one leaves both shaded blocks possible. We must distinguish
their union in the same nine tests.

\begin{center}
\begin{minipage}{.92\linewidth}\centering
\begin{tikzpicture}[x=1cm,y=1cm,font=\small]
\node at (2.07,2.45) {496 tested};
\node at (6.12,2.45) {477 untested};
\node[anchor=east] at (-.18,1.47) {4 tested};
\node[anchor=east] at (-.18,.37) {2 untested};
\fill[black!13] (4.14,.74) rectangle (8.1,2.22);
\fill[black!13] (0,0) rectangle (4.14,.74);
\draw (0,0) rectangle (8.1,2.22);
\draw (4.14,0)--(4.14,2.22); \draw (0,.74)--(8.1,.74);
\node at (2.07,1.48) {answer 2}; \node at (6.12,1.48) {answer 1};
\node at (2.07,.37) {answer 1}; \node at (6.12,.37) {answer 0};
\end{tikzpicture}
\captionof{figure}{The first test on six by 973. Each candidate pair occupies
one row and one column. The same 496 tested columns occur in two different answers.}
\label{fig:first}
\end{minipage}
\end{center}

Why test four of the six? Try five instead. Write \(x,y\) for the numbers
of tested and untested items in the other set. An answer of two leaves
the five paired with \(x\) items. An answer of one still allows those same
five to be paired with the other \(y\). Since \(n(9,5)=481\), both \(x\)
and \(y\) are at most 481, giving \(x+y\le962\). Splitting the six evenly
also falls short, for a different reason. After an answer of one, two
three-row blocks survive. Individually each could have \(n(9,3)=503\)
columns. The bound for searching them together is \(x+y\le966\).
Section~\ref{sec:six} derives the bounds for every possible first split.

Aigner \cite{A86} viewed a pair as an edge of a complete bipartite graph.
Its two vertex sets are its \emph{shores}, and a test counts the endpoints
of the edge in a chosen vertex set. Aigner determined the exact answers for
shores of sizes two through four. Li, Wu and Triesch \cite{LWT18} determined
the answer for five. Their first test selects three of the five when
eight tests are available and four when nine are available. With \(mn\)
possible pairs, any strategy needs at least \(\lceil\log_3(mn)\rceil\)
tests in the worst case. They proved that for each fixed \(m\) only finitely
many \(n\) attain this lower bound. Our first theorem gives the
exact searchable sizes for the next three shores.

\begin{theorem}[Exact searchable sizes]\label{thm:exact}
For six and seven items,
\[
\begin{aligned}
 n(k,6)&=2^k-k^2+7k-21 &&(k\ge9),\\
 n(k,7)&=2^k-\tfrac12(3k^2-25k+74) &&(k\ge14).
\end{aligned}
\]
At the six-item maximum every first test splits the six items four and two
for \(k\ge9\). For eight items, put \(B(k)=2^k-2k^2+16k-44\). Then
\[
 n(k,8)=
 \begin{cases}
 B(k),&11\le k\le12,\\
 B(k)+\lfloor(k-13)/3\rfloor,&k\ge13 .
 \end{cases}
\]
\end{theorem}

The seven-item polynomial gives 8,064 at thirteen tests, two above the true
maximum. For eight items the floor expression is one too small at twelve.
It agrees with \(B(k)\) from thirteen through fifteen. The first extra
item appears at sixteen, then a further item at nineteen. Supplement~S1
gives the complete initial values and all continuations for
\(1\le m\le8\), including the finite cases below the displayed ranges.

Our second result concerns the \emph{deficit} \(2^k-n(k,m)\). For six,
seven and eight items its leading terms are \(k^2\), \(3k^2/2\) and
\(2k^2\). The degree stays the same across these three sizes and rises at nine.

\begin{theorem}[Loss from the binary bound]\label{thm:leading}
For fixed \(m\ge2\) let \(d=\lceil\log_2m\rceil-1\). Then
\[
 \lim_{k\to\infty}\frac{2^k-n(k,m)}{\binom{k}{d}}=m-2^d .
\]
For \(m\ge3\) the error after subtracting \((m-2^d)\binom{k}{d}\)
from the deficit is \(O_m(k^{d-1})\). For \(m=2\) the deficit is one.
\end{theorem}

We prove the exact answers by constructing strategies and bounding what
every strategy can do. At eight items these two arguments explain different
features of the rounding. The construction locates the five groups of
candidate pairs that limit the gain. The upper bound explains why changing
the strategy cannot buy the next item sooner. The general leading term
has an analytic proof in Section~\ref{sec:general}.

\Needspace{8\baselineskip}
\section{How to certify a search}\label{sec:method}
\subsection{The three answers share the same pieces}
Write \(K_{h,w}\) for a block of candidate pairs with \(h\) rows and
\(w\) columns. A surviving state is a disjoint union of such blocks.
There is one unknown pair in the entire union. We call \(w\) the
\emph{width} and keep the first-shore rows fixed in this description even
if there are fewer columns than rows.

In block \(i\), choose \(b_i\) of the \(h_i\) rows and \(x_i\) of the
\(w_i\) columns. Write \(y_i=w_i-x_i\). The three states now possible are
\[
\begin{array}{c|l}
\text{answer}&\text{remaining pairs}\\ \hline
2&\bigsqcup_i K_{b_i,x_i}\\
1&\bigsqcup_i\bigl(K_{h_i-b_i,x_i}\sqcup K_{b_i,y_i}\bigr)\\
0&\bigsqcup_i K_{h_i-b_i,y_i}.
\end{array}
\tag{1}\label{eq:children}
\]
A block with zero rows disappears. A block with positive row count and
zero width represents no pairs and may be retained in sums. Each occurrence
of \(x_i\) comes from the same set of columns. Optimizing their sizes
separately for the different answers could prescribe two incompatible first
tests. We call the answer of one \emph{mixed} and the other answers
\emph{pure}.

\subsection{When the remaining rows are singletons}
A block with one row and \(w\) columns is a star with \(w\) leaves. A
union of these blocks is a star forest. We use two results about such
forests to finish strategies. Let \(N(r)\) be Aigner's searchable list of
star widths for \(r\) remaining tests, with entries \(g_r(i)\) in decreasing
order and zeros appended. The first entry is \(2^r\). For \(1\le j\le r\)
there are \(2^{j-1}\) entries \(\sum_{\ell=j}^r\binom r\ell\).
Deleting columns from distinct stars of this forest gives a strategy for
any list fitting in its slots \cite[Prop.~3.24]{A88}.

For a component with \(r\) tests left, define its width deficit as
\(d_i=2^r-w_i\). Set
\[
 L_u(r)=u2^r-\sum_{i=1}^{u}g_r(i)
       =\sum_{j\ge0}(u-2^j)_+\binom rj,\qquad (x)_+=\max(x,0).
 \tag{2}\label{eq:loss}
\]
Here \(L_0=L_1=0\) and \(\binom rj=0\) for \(j>r\). In a searchable star
forest every subset of \(u\) stars must have total deficit at least
\(L_u(r)\) \cite[Prop.~3.25]{A88}. In terms of widths, the \(u\) largest
sum to at most the sum of the first \(u\) entries of \(N(r)\). This is the
weak majorization condition. For a forest with at most eight stars,
checking all these inequalities also suffices for solvability, at every
budget \cite[Appendix A.2]{Kar26a}. This converse lets us stop specifying
a strategy as soon as all rows are singletons and the subset bounds hold.

A necessary bound for general blocks will also be useful. Give each row
of block \(i\) its own copy of the \(w_i\) columns. An original strategy
works on these disjoint stars by treating all copies of a column alike.
Consequently, if \(H=\sum_i h_i\), then
\[
 \sum_i h_i d_i\ \ge\ L_H(r).
 \tag{3}\label{eq:expanded}
\]

\subsection{What the certificates prove}\label{sec:certificates}
For a construction, fix the first few tests on the rows in each block.
Supply integer widths for their tested columns. Equation~\eqref{eq:children}
then fixes widths in all three children. Checking nonnegative widths at
every step and the subset bounds at every star terminal constructs a
strategy for every answer. The atlas in Supplement~S7 gives such choices.

For an upper bound, fix only the numbers of tested rows and leave the column
widths variable. Any surviving subfamily of blocks is itself searchable by
restriction. Its known upper bound is therefore a necessary inequality on
the very same \(x_i,y_i\) that appear in the other children. If a
nonnegative rational combination of child inequalities includes every
\(x_i,y_i\) with coefficient at least one, its right-hand side bounds
\(\sum_i(x_i+y_i)\). For example, the five-item argument in the introduction
adds two inequalities with coefficient one. An integer upper bound follows
by rounding the sum down.

An induction begins at zero tests, when at most one pair can remain. For
a fixed list of positive row counts, the largest sum of widths at zero
tests is one if that list contains a one and zero otherwise. At the next
budget the verifier checks a cover for every possible choice of tested
rows, identifying a test with its whole complement. When a child must be
expanded further, its new pieces still satisfy all ancestor equations.
Exact cancellation gives a contradiction at each excluded leaf. Thus an
upper certificate covers all possible tests, including those outside any
suggested construction. Supplement~S2 specifies the finite bases, symbolic
steps and the sources of the bounds below.

\section{Six and seven items}\label{sec:six}
\subsection{Why six forces a four-and-two split}
To see how a joint bound works, return to the balanced first test in
Section~1. It leaves two three-row blocks on an answer of one. In the
next test write \(x_i,y_i\) for the two column pieces of block \(i\).
Suppose we test no rows of the first block and two rows of the second.
An answer of zero still permits the three rows of the first block paired with
its \(y_1\) columns, so \(y_1\le n(8,3)=248\). An answer of one leaves
three rows paired with \(x_1\), one row with \(x_2\), and two rows with
\(y_2\). The eight-test induction gives their summed-width bound
\(x_1+x_2+y_2\le718\). Adding the two inequalities gives 966. Both
bounds constrain the same proposed test even though they come from
different answers. The other choices of rows for this second test have
covers of value at most 966 as well. Thus every balanced first test falls
at least seven columns short of 973.

For arbitrary budgets, write \(U_{\mathbf h}(r)\) for a certified upper
bound on the sum of widths of blocks with row counts \(\mathbf h\) and
\(r\) tests left. The induction of Supplement~S2 gives, for \(r\ge9\),
\[
\begin{aligned}
 U_{(4,2)}(r)&=2^{r+1}-r^2+5r-15, &
 U_{(3,3)}(r)&=2^{r+1}-r^2+4r-13,\\
 U_{(5)}(r)&=2^r-\tfrac12r^2+\tfrac32r-3, &
 U_{(6)}(r)&=2^r-r^2+7r-21 .
\end{aligned}
\]
Let \(F(k)=2^k-k^2+7k-21\). With \(r=k-1\), the mixed answer to a
four-and-two split gives \(U_{(4,2)}(r)=F(k)\). What can the other
splits achieve? If all six rows stay on one side, their two column pieces
each face a six-row problem, giving at most \(2U_{(6)}(r)\). If five stay
together, the same restriction gives \(2U_{(5)}(r)\). The balanced split
has the mixed bound \(U_{(3,3)}(r)\). Subtracting these bounds from \(F(k)\)
gives, respectively,
\[
 k^2-11k+37,\qquad 2k-11,\qquad k-3 .
\]
All are positive for \(k\ge10\). The nine-test base verifies the same
strict inequalities directly. Up to taking the complement of the entire
test, this exhausts the first splits.

To attain \(F(k)\), test four rows together with
\(n(k-1,4)=2^{k-1}-2(k-1)+2\) columns. The answer of two is an optimal
four-item problem. After four further tests on the longest specified
branch, every surviving block has one row. The five-test construction
in Supplement~S7 gives the column pieces at each step; the subset checks
in S3 certify every terminal for all \(k\ge9\). This proves both the
six-item formula and the necessity of its first split.

\subsection{Seven items and the sharp starting budget}
The simultaneous induction through eight rows supplies
\(n(k,7)\le 2^k-(3k^2-25k+74)/2\) from \(k=10\). The six initial tests
in the seven-row panel of Supplement~S7 attain this bound from fourteen
with the indicated continuations. A second construction lets us see the
root width directly from its first two tests. Put \(r=k-2\) and
\(Q=2^r\). Its four column pieces have widths
\[
 a=Q-1,\quad b=Q-3r+12,\quad
 c=Q-\tfrac12(r^2+3r-52),\quad w=Q-r^2+11r-55 .
\]
The first test uses two rows and \(b+c\) columns. On an answer of one
the state is \(K_{2,a+w}\sqcup K_{5,b+c}\). Test both rows and \(a\)
columns in the first block, and one row and \(c\) columns in the second.
A second answer of one leaves \(K_{2,w}\sqcup K_{1,b}\sqcup K_{4,c}\).
The root width \(a+b+c+w\) is the claimed polynomial. The certificate
completes this branch and both pure answers at each test.

At smaller budgets, the finite constructions and unrestricted exclusions
in Supplement~S1 give \(2^k-k^2+4k-10\) for \(4\le k\le10\) and
\(2^k-k^2+3k\) for \(11\le k\le13\). In particular the finite exclusion
at thirteen certifies that 8,063 is too large. This establishes the sharp
starting point independently of the choice of eventual construction.

\section{Eight items and the period of three}\label{sec:eight}
\subsection{Where the extra item fits}
Write \(B(k)=2^k-2k^2+16k-44\) as in Theorem~\ref{thm:exact}. For
\(k\ge15\), set \(g=\lfloor(k-13)/3\rfloor\). One fixed pattern of six
tests attains \(B(k)+g\) by varying the numbers of columns tested. For
example, its first test selects six of the eight rows together with
\(n(k-1,6)\) columns. Its answer of two is exactly the six-item problem
already solved. The full set of tests and widths is in Supplement~S7.

Follow six successive answers of one. There are eight stars left, one
per original row, and \(r=k-6\) tests remain. Figure~\ref{fig:eight}
isolates the part of the strategy that determines \(g\).
At fixed \(k\), increasing the gain parameter in the width expressions
by one reallocates three candidate columns from the eighth star to the
second. The total width of the eight stars stays constant. The widths
of stars 2 through 6 grow by three in total. Their deficits, in row order, are
\[
 \tbinom k2-6k+11-3g,\quad k-5,\quad k-5,\quad k-3,\quad k-4 .
\]
By \eqref{eq:loss} their sum must be at least \(L_5(r)=4+3r+\binom r2\).
Subtract this minimum. Exactly \(k-13-3g\) units remain. Each further item
at the root would spend three of them in this assignment.

\begin{center}
\begin{minipage}{.94\linewidth}\centering
\begin{tikzpicture}[x=1.02cm,y=.75cm,font=\small]
\fill[teal!16] (1,0) rectangle (6,1);
\draw[step=1] (0,0) grid (8,1);
\foreach \i in {1,...,8} \node at (\i-.5,1.33) {\i};
\node[anchor=east] at (-.15,1.33) {row};
\node[anchor=east] at (-.15,.5) {width change};
\foreach \i/\v in \eightwidthchanges
  \node at (\i-.5,.5) {\(\v\)};
\draw[decorate,decoration={brace,mirror,amplitude=4pt}] (1,-.12)--(6,-.12);
\node at (3.5,-.68) {five widths gain three in total};
\end{tikzpicture}
\captionof{figure}{Width changes per unit of \(g\) after six answers of
one. The shading selects the same five stars as the hatching in the full
eight-item diagram in Supplement~S7.}\label{fig:eight}
\end{minipage}
\end{center}

The total deficit is \(4\binom k2-18k+55=L_8(k-6)\), independent of \(g\).
Thus the inequality for all eight stars would miss the rounding. For this
pattern a rational combination of terminal inequalities gives
\(3D\ge6k^2-49k+145\) for the root deficit \(D\); the five-star inequality
has coefficient \(1/3\). The assignment attains the integer ceiling. Its
remaining piece capacities and \emph{all} terminal subset inequalities
are checked on each residue class of \(k\) as explained in Supplement~S3.
Other assignments can place slack at a different terminal. At eleven through
fourteen tests separate constructions attain \(B(k)\).

\subsection{Why no other test can gain more}
Suppose \(k\ge11\) and a strategy can search \(B(k)+h\) columns with integer
\(h\ge1\). Every first test capable of a gain must split the eight rows
two and six. Suppose its mixed answer leaves \(K_{6,x}\sqcup K_{2,y}\).
One pure answer contains \(K_{6,y}\), so \(y\le n(k-1,6)\). The next test
must keep the two rows together and split the six into four and two.
Supplement~S2 derives these reductions from bounds on all the other cuts.

After two mixed answers, three blocks have four, two and two rows. Mark
the two-row block that shares its parent with the four-row block. The two
blocks of equal height play different roles: the restriction on \(y\)
comes from the \emph{other} block's ancestor. Put \(t=k-2\) and let \(w\)
be the piece of \(y\) sent to a pure two-row component by the second test.
The unused widths at the two pure answers are
\[
 \delta=n(k-1,6)-y\ge0,\qquad v=2^t-1-w\ge0 .
\]
The marked induction of Supplement~S2 establishes the following three
facts. At entry its integer parameter is
\(\Lambda=t-8-3h-\delta-2v\). Every return to the marked three-block state
reduces \(\Lambda\) by at least ten. Every successful exit requires
\(\Lambda\ge3\). A finite base excludes survival to small remaining
budgets, so every successful branch eventually exits.

One unit in the exit condition matters. In a pure sibling there are three
rows opposite \(2^s-b\) columns and two rows opposite \(2^s-a\).
The weighted cover gives \(2a+b\ge\binom{s+1}{2}+3\). For \(s\ge5\),
equality forces \(a=1\). If a proposed exit has \(\Lambda\le2\), its
shared widths instead force \(a\ge2\). Applying the strengthened
inequality rules out this last case. Supplement~S2 gives the equality
argument and the induction covering every exit.

After \(j\) returns, the entry and exit bounds give
\[
 3h+\delta+2v+10j\ \le\ t-11=k-13 .
 \tag{4}\label{eq:gain}
\]
In particular \(3h\le k-13\) for \emph{any} strategy with positive gain.
At eleven and twelve this excludes every gain. From thirteen it gives
the floor in Theorem~\ref{thm:exact}. The constructions attain the bound
in each range, completing the proof.

\section{How much is lost for an arbitrary fixed shore?}\label{sec:general}
The exact calculations distinguish strategies that differ by a single
column. For the leading term it suffices to obtain two bounds whose
highest-degree terms agree. The upper bound must follow the several blocks
a test leaves behind, even when the search starts with just one.

Fix a list of positive row counts \(\mathbf h=(h_1,\ldots,h_s)\), write
\(H=\sum_i h_i\), and let \(W_{\mathbf h}(k)\) be the largest sum of
widths searchable in \(k\) tests, allowing zero widths. There are \(s\)
individual binary bounds. Their combined deficit for a given assignment
is \(D=s2^k-\sum_i w_i\).

\begin{lemma}\label{lem:uniform}
Put \(C=H-s\). For \(H\ge1\) and \(k\ge C\),
\[
 W_{\mathbf h}(k)\le s2^k-L_H(k-C).
\]
\end{lemma}
\begin{proof}
The shift \(C\) counts the rows in excess of one per block. Whenever a
test splits a block's rows, both pieces appear on the mixed answer and
that shift decreases. A second identity treats tests leaving every block's
rows together. By Pascal's identity applied to \eqref{eq:loss},
\[
 L_H(t)-L_H(t-1)=
 L_{\lfloor H/2\rfloor}(t-1)+L_{\lceil H/2\rceil}(t-1).
\]
The function \(u\mapsto L_u(t-1)\) is convex on integers, being a sum of
positive parts. Hence for any \(a+b=H\) the right-hand side is at most
\(L_a(t-1)+L_b(t-1)\).

We induct on \(k\) simultaneously for all row lists, starting each at its
own \(C\). For \(H=1\) the assertion is the binary bound. For \(H\ge2\),
each \(h_i\le C+1\). At \(k=C\), \eqref{eq:expanded} gives
\[
 D\ge\frac{L_H(C)}{C+1}
 \ge\frac{(H-1)+(H-2)C}{C+1}>H-2 .
\]
Since \(D\) is an integer, \(D\ge H-1=L_H(0)\).

For the induction step keep the entire mixed answer and the pure pieces
of the blocks whose rows were left together. These subfamilies contain
every \(x_i,y_i\) once in \(2s\) formal blocks, so their deficits sum to
\(D\). If \(I\ge1\) blocks have their rows split, the mixed child has
\(H\) rows, \(s+I\) blocks and shift \(C-I\). Induction there gives at
least \(L_H(k-1-C+I)\ge L_H(k-C)\); the other deficits are nonnegative.

If no rows split, the mixed child keeps the original profile. The two
pure subfamilies partition the rows into counts \(a+b=H\), each with
shift at most \(C\). Writing \(t=k-C\ge1\), their three deficits sum to
at least \(L_H(t-1)+L_a(t-1)+L_b(t-1)\ge L_H(t)\). An empty subfamily
contributes zero. This proves the lemma.
\end{proof}

\subsection{A strategy that attains the leading term}
How can we build nearly \(2^k\) columns per block while finding which row
contains the first defective? Give all \(H\ge2\) rows distinct binary
labels of length \(q=\lceil\log_2H\rceil\). Divide the columns in each
block into \(2^q\) groups, one per binary label. In test \(j\), include
every row and column whose \(j\)-th bit is one. Once the \(q\) sums have
been reported, any surviving row determines its partner's column label
bit by bit. Distinct surviving rows in the same block require distinct
labels. The state is now a star forest.

Set \(d=q-1,\ c=H-2^d,\ t=k-q\) and take \(k\ge2q\). Assign the rows
distinct ranks \(1,\ldots,H\). In each row's block give the column group
with its complementary label exactly \(g_t(i)\) columns, where \(i\) is
the row's rank. Give each remaining group \(g_t(2^d)\) columns. If all
answers are one, the stars are precisely the first \(H\) slots of \(N(t)\).
Any other answer fixes at least one bit of the row, leaving at most \(2^d\)
rows. For each \(j<d\), at most \(2^j\) surviving widths can exceed
\(g_t(2^j+1)\): only widths assigned to ranks \(1,\ldots,2^j\) can do so.
The step lengths of \(N(t)\) therefore let all surviving stars fit in
distinct slots. Aigner's construction finishes in every case.

Summing the widths of the \(H\) designated groups and the \(s2^q-H\)
remaining ones gives
\[
 s2^k-\sum_{j=0}^{d-1}(s2^q-2^j)\binom{k-q}{j}
       -c\binom{k-q}{d}.
\]
Equation~\eqref{eq:loss} shows that the lower bound \(L_H(k-C)\) on the
deficit has the same leading coefficient \(c\) in the binomial basis.
This proves Theorem~\ref{thm:leading} for one block, and more generally
\[
 s2^k-W_{\mathbf h}(k)
  =(H-2^d)\binom{k}{d}+O_{\mathbf h}(k^{d-1})\qquad(H\ge3).
\]
For \(H=1\) binary search gives zero loss. For \(H=2\) the loss is one
from \(k=1\): the construction covers \(k\ge2\), and at one test
\(K_{2,1}\) and \(K_{1,2}\sqcup K_{1,1}\) realize the two possible
profiles. The lemma also gives the finite bound
\(n(k,m)\le2^k-L_m(k-m+1)\) for \(k\ge m-1\). Supplement~S4 gives the
corresponding leading term when the value per column varies from block to
block.

\subsection{What remains open}
The binary labels explain the degree of the loss for every fixed size.
They leave open the lower-order terms for nine or more rows. At eight,
one fixed pattern of initial row tests suffices from fifteen tests onward.
For any prescribed finite family of such patterns the optimal deficit is
eventually a polynomial plus a bounded periodic term (Supplement~S5).
Does some finite family attain the \emph{unrestricted} optimum for every
sufficiently large budget, for each fixed \(m\)? An affirmative answer would
extend this description beyond the eight-item case.

There is also a way to see the scale of the benefit of adaptivity. For
fixed \(m\ge3\), requiring every test to be chosen in advance changes the
loss from the binary bound from polynomial to exponential in \(k\). Both
searchable sizes still have ratio tending to one relative to \(2^k\).
Supplement~S6 proves the comparison using binary codes and gives the
equivalent search interpretations.

\paragraph{Data and code.}
The constructions, exclusions, finite bases and exact verifiers are
archived on Zenodo (DOI \href{https://doi.org/10.5281/zenodo.23996093}%
{\mbox{10.5281/zenodo.23996093}}). The archive maps each claim to its
proof. Sections S1--S3 of the supplement identify the inputs for each
range. Earlier finite results are in \cite{Kar26b}; the star-forest
continuation theorem is in the revised companion \cite{Kar26a}.

\paragraph{Acknowledgments.} The author thanks Konstantin Knop and Igor
Krivokon for introducing the problem and for their advice and consultations
over the years. This research received no external funding. The author
declares no competing interests.

\paragraph{Use of generative AI.} The author used OpenAI Codex and, to a
lesser extent, Anthropic Claude for literature search, development and
exposition of mathematical arguments, software development and debugging,
computational experiments and manuscript preparation. Codex also assisted
with diagram code, layout and revision of the exposition. The author
directed the work and assumes responsibility for the claims, references,
software, computations and figures. Various model versions were used; a
complete historical version record is unavailable.

\bibliographystyle{plain}
\bibliography{fixed-shore}
\end{document}

\thispagestyle{plain}
\begin{center}
{\small Supplement to\par}\vspace{3pt}{\Large\papername\par}\vspace{.6em}
\papersignature
\end{center}

\section{Exact values and their proof sources}\label{sec:table}
For one defective in each of two disjoint sets, \(n(k,m)\) is the largest
second set searchable with \(k\) adaptive count tests when the first set
has \(m\) items. Table~\ref{table:1} gives the complete beginning of each
sequence. For \(m\le5\) the continuations of Aigner and of Li, Wu and
Triesch \cite{A86,LWT18} are
\[
\begin{aligned}
 n(k,1)&=2^k &(k&\ge0), & n(k,2)&=2^k-1 &(k&\ge0),\\
 n(k,3)&=2^k-k &(k&\ge2), & n(k,4)&=2^k-2k+2 &(k&\ge3),\\
 n(k,5)&=2^k-\tfrac12(k^2-3k+6) &(k&\ge11).
\end{aligned}
\]
At six and seven the continuations are \(2^k-k^2+7k-21\) from nine and
\(2^k-(3k^2-25k+74)/2\) from fourteen, respectively. At eight the
continuation is \(2^k-2k^2+\lfloor(49k-145)/3\rfloor\) from thirteen.

\input{initial-values}

\paragraph{Finite bases.} The compact finite slice in the
\href{https://doi.org/10.5281/zenodo.23996093}{fixed-shore evidence archive}
replays the claims used from the earlier finite proof \cite{Kar26b} through
eight tests. It includes positive strategies and negative bounds with
their required lower-budget premises. If the displayed first shore is
larger, transposition and restriction transfer the corresponding bounds.
At nine and ten tests the eight-row positive witnesses establish widths
447 and 939; separate unrestricted exclusion certificates reject 448
and 940. For seven rows the finite package verifies the earlier maxima
at nine through thirteen. In particular it supplies both a strategy at
8,062 and a refutation at 8,063 for thirteen tests. The later six- and
eight-row ranges use the proofs in \S\S\ref{sec:covers}--\ref{sec:positive}.

The archive's claim-to-proof map identifies the minimal replay for each
assertion. Its group names U6 and U8 denote the simultaneous width
inductions, U8-marked the unrestricted marked-state exclusion, F7 the
finite seven-row campaign, B8 the positive eight-row boundary, and S, P
and F8 the symbolic constructions and their checks. These names allow
the proofs below to be located in the archived data independently of
section numbering in the article.

\section{Unrestricted upper certificates}\label{sec:covers}
\subsection{Covers and extension to every budget}
For each block split its nonnegative width into pieces \(x_i,y_i\).
A child with \(r-1\) tests left obeys every previously proved inequality
on its surviving subfamilies. Write these inequalities in the same piece
variables \(z_v\) as \(\sum_v c_{qv}z_v\le U_q(r-1)\). A certificate gives
nonnegative rational multipliers \(\lambda_q\) with
\(\sum_q\lambda_qc_{qv}\ge1\) for each \(v\). Then
\[
  \sum_i(x_i+y_i)\le
   \left\lfloor\sum_q\lambda_q U_q(r-1)\right\rfloor .
\]
For a weighted objective the required coefficients on the left are its
specified weights. Each cited bound has its own ordered row counts and
validity range. The zero-test base is one if the list contains a one-row
block and zero otherwise. All row cuts must be covered. Complementing the
whole test exchanges answers zero and two and preserves one. Independently
complementing the test in one block can change the surviving union.

The finite induction U6 ends at budget nine. Its symbolic continuation
uses expressions \(a2^r+br^2+cr+d\). To certify a nonnegative difference
from \(r=10\) onward, the checker verifies its value and its first two
forward differences at ten are nonnegative and \(a\ge0\). Its third
forward difference is \(a2^r\), so induction in \(r\) proves all later
inequalities. U8 supplies the auxiliary bounds for profiles through
eight rows and the seven-item polynomial upper bound from ten tests
onward. For a deeper refutation the checker reconstructs the three
children and retains ancestor equations. At a leaf it checks an exact
positive contradiction using nonnegative multipliers on inequalities
and unrestricted multipliers on equalities. At an expanded node it
verifies that every row-cut orbit is present. A finite positive strategy
and a failed attempt to extend that strategy could not supply this coverage.

\subsection{The marked induction for eight rows}
Let \(B(k)=2^k-2k^2+16k-44\), \(k\ge11\), and suppose width \(B(k)+h\)
is searchable for some integer \(h\ge1\). The U8 profile covers exclude
each first cut except two plus six. Write its mixed child as
\(K_{6,x}\sqcup K_{2,y}\). The pure answers are \(K_{6,y}\) and \(K_{2,x}\),
so with \(r=k-1\) and \(\delta=n(r,6)-y\ge0\),
\[
 x=B(k)+h-n(r,6)+\delta,\qquad y=n(r,6)-\delta .
\]
At the second test the individual five- and six-row bounds exclude the
\(1+5\) and \(0+6\) splits of the six. The mixed bounds for row profiles
\((4,2,1,1)\) and \((3,3,1,1)\) exclude splitting the two rows internally.
Adding the \((3,3,2)\) mixed bound and the pure two-row bound excludes
\(3+3\). Both remaining orientations of the \(4+2\) split are retained.

Put \(t=k-2\), and let \(w\) be the piece of \(y\) in the pure two-row
component at this second test. Restriction gives \(v=2^t-1-w\ge0\).
In the mixed state mark the two-row component descended from the six
rows. The width deficits of this block and its four-row sibling
are denoted by \(D_2,D_4\); the third is \(D'_2\). The marked pair's
widths sum to \(x\), and all three widths sum to \(x+y-w\). Thus
\[
 D_4+D_2=2^{t+1}-x,\qquad \sum_iD_i=3\cdot2^t-(x+y-w).
\]
For the row profiles \(F=(4,2,2)\) and \(H=(4,2,1,1)\) define
\[
 S_F(t)=2t^2-9t+22,\qquad S_H(t)=2t^2-7t+17,\qquad
 z=\sum_iD_i-S_P(t)\quad(P=F,H).
\]
In either profile put
\(E=t^2-3t+5-(D_4+D_2)\), keeping the same marked ancestry. The region
under consideration is \(E\ge1,\ 0\le z\le t-4,\ t\ge4\). The profile
width bounds give \(z\ge0\) in any searchable state in this budget range.
At entry the two width identities give
\[
 E=h+\delta,\quad z=t-3-h-v,\quad
 \Lambda=2z-E-t-2=t-8-3h-\delta-2v.
\]
In particular the state enters the \(F\) region.

For \(H\), all tests except the row choice \((0,2,1,1)\) have
contradiction covers. The survivor returns to \(H\) with one fewer test,
\(z'\le z-4\) and \(E'\ge E+1\). The base at four tests is impossible.
At five through seven the return would force \(z'<0\); at larger budgets
a searchable return has \(0\le z'\le t-8\le(t-1)-4\).
Induction excludes the whole \(H\) region.

The base for \(F\) at four is also impossible. At \(t\ge5\) only
three row choices remain, in four/marked-two/other-two order. Choice
\((0,2,1)\) enters \(H\) with \(E'\ge E+1,\ z'\le z-1\), inside the
excluded region if searchable. Choice \((0,2,2)\) returns to \(F\)
with \(E'\ge E+1,\ z'\le z-5\). Again a searchable return has
\(0\le z'\le (t-1)-4\). Substitution gives \(\Lambda'\le\Lambda-10\).
The third choice, \((1,2,0)\), is the exit. Its mixed answer has profile
\((3,2,2,1)\).

\paragraph{The last exit.} A cover adds four inequalities: the pure
bound for the marked two, a mixed \((3,2)\) bound, the full mixed bound,
and the opposite pure weighted \((3,2)\) bound. Their multipliers are
\(3/2,1/2,1,1/2\), respectively. The result
\[
 \sum_iD_i+\tfrac12(D_4+D_2)\ge\tfrac12(5t^2-20t+51)
\]
implies \(\Lambda\ge0\). The next exclusion eliminates \(\Lambda\le1\)
and closes all \(F\) budgets from four through twelve. For a later exit
with \(\Lambda\le2\), assuming the two-row deficit of the pure \((3,2)\)
sibling is at most one gives a rational contradiction with positive
margin \(1/2+t/6\). Thus that deficit is at least two.

To see the improvement this permits, write the sibling as
\(K_{3,2^s-b}\sqcup K_{2,2^s-a}\). For \(s\ge1\) a weighted induction
gives \(2a+b\ge\binom{s+1}{2}+3\). At \(s\ge5\), five of its six
cut orbits give at least one further unit. For the sixth, \((1,1)\),
the mixed profile \((2,1,1,1)\) implies \(a+b\ge\binom{s+1}{2}+2\).
Since \(a\ge1\), equality in the weighted bound forces \(a=1\).
Here \(a\ge2\), so the right-hand side improves to \(\binom{s+1}{2}+4\).
With this bound the final certificate covers all cuts of the exit's mixed
child and every further expanded branch, excluding \(\Lambda\le2\).
Each symbolic exclusion has finite base cases and exact tail inequalities.

Consequently the last exit needs \(\Lambda\ge3\). Returns reduce the
remaining budget, so they cannot continue forever. Subtracting ten per
return from the entry value yields \(3h+\delta+2v+10j\le k-13\) for
\(j\) returns. In the archive, U8-marked replays the underlying profile,
weighted, marked-region and sibling-equality certificates together. This
last-exit condition accounts for the constant inside the floor.

\section{Construction certificates and alternative strategies}\label{sec:positive}
\subsection{Checking one assignment for infinitely many budgets}
A prefix specifies the initial row cuts, including all three children at
each node, until only stars remain. A component of deficit \(d\) at \(r\)
tests is split into widths \(2^{r-1}-a,\ 2^{r-1}-c\), with \(a+c=d\).
All descendants using the same piece inherit its value. A certificate
supplies the tested deficits \(a\) and the root deficit. The checker
recovers \(c\), validates both piece capacities and checks every subset
of each terminal using \(L_u(r)=\sum_{j\ge0}(u-2^j)_+\binom rj\).
For at most eight stars this last check suffices by \cite[Appendix A.2]{Kar26a}.

For a base \(b\) and period \(p\), write \(k=b+pu\). All deficits are
given in the binomial basis \(1,u,\binom u2,\ldots\). Shared-piece
identities hold coefficient by coefficient. Nonnegative binomial
coefficients of each slack certify it for every \(u\ge0\). To check
binary capacities against a polynomial, expand
\(2^{r+pu}=2^r\sum_{j\ge0}(2^p-1)^j\binom uj\); only finitely many
coefficients can be affected by subtracting the polynomial. These
are sufficient tests for the stated assignments. The verifiers use exact
arithmetic for all identities, covers and capacities.

Group S checks prefixes of depths zero through three for one through four
rows, valid from budgets zero through three respectively. The five-, six-
and seven-row prefixes have depths four, five and six and start at eleven,
nine and fourteen. Group P adds a strict eight-row baseline from twelve
and the seven-row construction at thirteen. The eight-row single-floor
assignment F8 has root deficit
\(44-14k+4\binom k2-g\) with \(g=\lfloor(k-13)/3\rfloor\).
Testing its three residue classes from fifteen proves all capacities,
all terminal subset inequalities and the rational dual combination in the
article. Substituting fourteen into this assignment violates a terminal
subset inequality; B8 supplies separate ordinary positive strategies at
eleven through fourteen. The restriction is on this particular assignment.

\subsection{Other ways to attain the same widths}
For six rows one can also finish the outer prefix by recursion. Write
\(A=2^r\). Two searchable labelled star lists are
\[
\begin{aligned}
 T_r&=(A-2,A-2,A-r+2),\\
 P_r&=(A-r-1,A-r-2,A-r-2,A-2r+3,\\[-3pt]
    &\hspace{5em} A-(r^2-5r+8)/2,A-(r^2+r-2)/2).
\end{aligned}
\]
Use the certified bases at \(r=2\). For \(r>2\), retain \(T_{r-1}\) on
the mixed answer. The pure pieces \(2^{r-1},2^{r-1},2^{r-1}-1\) split
into a single star and a canonical pair. For \(P_r\) retain \(P_{r-1}\)
with its first two labels interchanged. The six pure pieces, with
\(Q=2^{r-1}\), are \(Q,Q-2,Q-1,Q-2,Q-r+3,Q-r\). Labels \(1,3,6\)
fit in distinct canonical slots and the others form \(T_{r-1}\).
The direct outer prefix uses \(P_r\) for its only nonembedded continuation
from eleven tests; ordinary trees supply nine and ten. The S construction
in the atlas covers all \(k\ge9\) in a single family.

For the explicit seven-row construction in the article, put \(u=k-14\).
Its root deficit is \(3\binom u2+31u+156\). The remaining assignments
are in the same binomial basis; every final state embeds in canonical
slots or has a known small-shore continuation. This gives a second
construction from fourteen alongside the one in the atlas.

An alternative eight-row construction iterates a three-test return.
Write \(k=19+3u+t\) with \(u,t\ge0\). Four canonical widths at budget
\(t\) are \(A=2^t,\ B=A-1,\ C=B-t,\ D=C-\binom t2\).
At the root their multiplicities satisfy
\[
 D_u=4,\quad C_u=58+12u,\quad B_u=398+167u+18u^2,\quad
 A_u=2^{19+3u}-B_u-C_u-D_u .
\]
The resulting root width is \(B(k)+2+u\), where \(B(k)\) is the eight-row
baseline of the article. The certificate assigns compatible multiplicities
at every internal component. Two star terminals have continuations that
return to parameter \(u-1\) in three tests, with explicit pure siblings;
the remaining terminals have direct small-star continuations. Taking
\(u=\lfloor(k-19)/3\rfloor\) attains the formula from nineteen. The
archive contains the starting matrices, earlier constructions and exact
return checks. This is an alternative lower proof; S2 supplies the upper one.

\section{The leading deficit for weighted widths}\label{sec:weighted}
Take blocks with row counts \(h_1,\ldots,h_s\), total \(H\ge2\), and
positive real weights \(\lambda_i\). Let \(W_{\mathbf h,\lambda}(k)\)
be the largest searchable value of \(\sum_i\lambda_i w_i\).
Write \(d=\lceil\log_2 H\rceil-1,\ c=H-2^d\). Allocate \(c\) units
among the blocks with capacities \(h_i\), filling the least expensive
blocks first. Let \(c_i^*\) be the allocation and \(B=\sum_i\lambda_ic_i^*\).

\begin{theorem}
The weighted leading deficit is
\[
 \lim_{k\to\infty}
 \frac{(\sum_i\lambda_i)2^k-W_{\mathbf h,\lambda}(k)}{\binom{k}{d}}=B .
\]
\end{theorem}
\begin{proof}
Order the weights increasingly and write \(a_i=2^k-w_i\).
Summation by parts gives
\[
 \sum_i\lambda_i a_i=\lambda_1\sum_i a_i+
    \sum_{j=2}^s(\lambda_j-\lambda_{j-1})\sum_{i=j}^s a_i .
\]
Each tail is a restriction to \(H_j=\sum_{i=j}^s h_i\) rows. By the
uniform width lemma in the article its deficit has lower bound with
degree-\(d\) coefficient \((H_j-2^d)_+\) in the binomial basis.
Consequently the necessary leading coefficient is at least
\(\lambda_1 c+\sum_{j=2}^s(\lambda_j-\lambda_{j-1})(H_j-2^d)_+=B\):
the excess that cannot fit in cheaper blocks must be paid for in this tail.
In the binary-label construction assign \(c_i^*\) ranks above \(2^d\)
to block \(i\). Its leading loss is then \(c_i^*\binom{k}{d}\), and all
other losses have smaller degree. At \(d=0\) there are no smaller-degree
losses. This gives the matching upper bound for the deficit.
\end{proof}

\section{What a finite family of test patterns can attain}\label{sec:templates}
A \emph{prefix template} specifies the row cuts at the nodes of a finite
decision tree starting from one \(m\)-row block. Every terminal is a star
forest. The tested column widths are unknown integers shared across the
three children. At remaining budget \(r\) impose the piece capacities and
the inequalities \(\sum_{i\in I}d_i\ge L_{|I|}(r)\) on all star subsets.
The depth \(q\) counts tests along the longest specified branch.

\begin{proposition}\label{prop:equivalence}
Suppose \(k\ge q\) and star majorization suffices at every terminal's star
count and remaining budget. A width assignment is realizable with that
prefix exactly when it satisfies the shared integer system.
\end{proposition}
\begin{proof}
Every strategy gives the same piece size to all occurrences of a piece,
and satisfies each necessary subset bound. Conversely, an integer
assignment selects a legal test at each node. Its children inherit the
same selections, and the assumed sufficiency finishes each terminal.
\end{proof}
There are at most \(m\) terminal stars on a branch. The cited converse
applies at every budget for at most 29 stars, and for arbitrary support
when the remaining budget is at most five. It also applies for at most
31, 32, or 33 stars with at least seven, eight, or nine tests left,
respectively. The result through sixteen stars has an analytic proof.
For arbitrary fixed \(m\), \(k-q\ge(m+1)^{6m}\) suffices
\cite[Theorems 1.1--1.3 and Appendix A.2]{Kar26a}.

Call a cut \emph{strict} if it splits the rows of at least one block.
With \(H\) rows in \(s\) blocks the integer \(H-s\) decreases on all
three answers of a strict cut. On the mixed answer each internal split
adds a block while preserving its row total. On a pure answer an internally
split block loses at least one unit of the potential and no other block
gains. There are finitely many strict templates from \(m\) rows, each of
depth at most \(m-1\). For \(m\le8\) their best feasible root width is
\(n(k,m)\) at every \(k\ge13\). Indeed S and P in \S\ref{sec:positive}
give attaining strict prefixes, including seven rows at thirteen and eight
rows at thirteen and fourteen. The classical frontiers through five and
the unrestricted upper bounds in the article close the reverse inequality.

\begin{theorem}\label{thm:periodic}
For any fixed nonempty finite family \(\mathcal T\) of prefix templates
let \(F_{\mathcal T}(k)\) be its least root deficit under the terminal
majorization conditions. Eventually \(F_{\mathcal T}(k)=P(k)+\eta(k)\),
where \(P\) is a rational polynomial and \(\eta\) is bounded and periodic.
The resulting assignments yield actual strategies for all sufficiently
large \(k\).
\end{theorem}
\begin{proof}
Fix a template \(T\) of depth at most \(q\). At depth \(a\) assign every
component deficit \(2^{q-a}L_m(k)\), splitting it equally at every test.
Every terminal subset condition holds since each component deficit is at
least \(L_m(k)\). The root cost \(C(k)=2^qL_m(k)\) is polynomial in \(k\)
and eventually at most \(2^{k-q}\), so every capacity holds. For \(m=1\)
the cost is zero. In a nonnegative split \(a+c=d\), each descendant
deficit is at most the root deficit. We may therefore bound all variables
by \(C(k)\) and drop the exponential caps at sufficiently large \(k\).
The remaining system has a constant matrix and polynomial right-hand
sides. Its rational optimum \(P_T(k)\) is eventually a single polynomial:
there are finitely many bases, their coordinates are polynomial, and their
feasibility signs and cost comparisons stabilize.

To bound the cost of integrality, prescribe nonnegative \emph{integer}
minimum deficits at the terminal components. Work up the tree. For an
internal component let \(A\) be the maximum of the requirements on
its tested piece, taken over all child occurrences, and \(C\) the
corresponding maximum on its untested piece. Its requirement is \(A+C\).
If every computed requirement fits its component capacity, every integer
root deficit between its requirement and capacity can be allocated.
Indeed, at a component with capacity \(2Q\) and \(A+C\le d\le2Q\),
set \(a=\max(A,d-Q)\) and \(c=d-a\). Then \(A\le a\le Q\),
\(C\le c\le Q\). Copy each chosen piece to its child occurrences and
continue. This proves the integer reconstruction used below.

Round the terminal components of a rational optimum upward. Each initial
error is in \([0,1)\). At an internal component the error in the
reconstructed minimum is the maximum child error over uses of its first
piece plus the maximum over uses of its second, because the rational
values of repeated pieces coincide. Starting instead with one at each
terminal and using the same recursion gives a bound \(R_T\le2^q\) at
the root. Computed requirements remain polynomially bounded, so they fit
eventually. Rounding up preserves every terminal subset inequality.
If the resulting root cost exceeds the auxiliary cutoff \(C(k)\), use
the integer assignment of cost \(C(k)\) above.
Writing \(F_T\) for the integer optimum gives
\[
 \lceil P_T(k)\rceil\le F_T(k)\le\lceil P_T(k)\rceil+R_T-1 .
\]

After clearing fixed denominators, these polynomially bounded systems
have eventually quasipolynomial optima and optimal assignments by
\cite[Corollary 1.16, Property 3a]{BGW17}; see also \cite{Sh18}.
A quasipolynomial is polynomial on each residue class of some fixed
modulus. Classical proximity for a fixed integral matrix also gives a
bounded gap \cite[Corollary 2]{CGST86}; the reconstruction above gives
the tree-specific bound. The smallest rational optimum in the finite
family is eventually one polynomial \(P\). The bounded gap for a rational
minimizer bounds \(F_{\mathcal T}-P\) from above, and relaxation bounds
it below by zero. On each common residue class this difference is
eventually polynomial and bounded, hence constant. Finally the fixed-star
eventual converse and Proposition~\ref{prop:equivalence} realize the
assignments.
\end{proof}

For \(m>1\) the degree is at most \(\lceil\log_2m\rceil-1\) by the
polynomial budget \(C(k)\). The terminal budget for realizability does
not determine the onset or period of the quasipolynomial. Whether some
finite family attains the unrestricted optimum for each fixed \(m>8\)
remains open.

A smaller class of prescriptions already fails at eight. A fixed entry of
\(N(k-q)\) with fixed offset \(q\) is \(2^{k-q}-\sum_{j<\ell}\binom{k-q}{j}\)
eventually for some fixed \(\ell\). A sum with fixed nonnegative integer
multiplicities has form \(c2^k-P(k)\). For an eventual attaining width,
\(c\) must equal one. Of finitely many such polynomials one eventually
dominates. But for \(k\ge13\) the eight-item deficit is
\[
 (6k^2-49k+145)/3+\bigl((k-13)\bmod3\bigr)/3.
\]
No finite list of those fixed sums attains it for all large \(k\).
The six-test pattern in the article allows its widths to vary with the
residue class.

\section{Coding interpretations and tests fixed in advance}\label{sec:coding}
In a binary adder channel with complete feedback, each sender sees its own
message and the previous outputs and sends one bit. The output is the sum
of the bits. Identify a message with an item in its shore. The encoding
functions specify exactly the two tested subsets at each history.
Conversely every test supplies these encoding functions. Distinct
transcripts for every allowed pair are equivalent to zero-error recovery;
an early finish can be padded with empty tests. This gives the
unequal-message version of the correspondence of Florin, Ho and Jiang
\cite{FHJ22}. A disjoint union of blocks corresponds to a promise on
the allowed pairs. For binary inputs the unordered set of the two symbols
also determines their sum and conversely. Florin, Ho and Jiang determine
the symmetric asymptotic rate. With fixed \(m\) and growing \(k\), the rate
pair approaches \((0,1)\), and the additive deficit distinguishes strategies
with the same limiting rate.

There is an exact interpretation in binary black-peg Mastermind
\cite{MS24}. For a known weight-two secret \(z\in\{0,1\}^N\) and test
vector \(x\), put \(y=x\cdot z\). The black-peg score is
\(N-|x|-2+2y\) and the Hamming distance is \(|x|+2-2y\).
Thus each query and answer recovers the other model's answer, with the same
allowed pairs, adaptivity and worst-case depth.

Let \(n_{\rm NA}(k,m)\) denote the largest opposite shore searchable by
\(k\) tests fixed in advance and set \(D_{\rm NA}=2^k-n_{\rm NA}(k,m)\).
In this case both shores use fixed binary codebooks \(A,B\subseteq\{0,1\}^k\);
the integer sums \(a+b\) must all be distinct. Write \(d_*(k,m)\) for the
largest possible minimum distance of \(m\) binary words of length \(k\).

\begin{proposition}
For \(2\le m\le2^k\),
\[
  2^{k-d_*(k,m)}\le D_{\rm NA}\le \binom m2\,2^{k-d_*(k,m)}.
\]
For fixed \(m\ge2\), \(D_{\rm NA}=\Theta_m(2^{\alpha_m k})\), where
\(\alpha_m=1-\lfloor m^2/4\rfloor/\binom m2\).
\end{proposition}
\begin{proof}
For fixed \(A\), join \(b,b'\) when \(a+b=a'+b'\) for distinct \(a,a'\in A\).
The allowed \(B\)'s are precisely the independent sets of this conflict
graph \cite[Theorems 1.3.5 and 1.4.1]{C83}; see \cite{T83} for prescribed
codebooks. Two words at distance \(h\) give \(2^{k-h}\) disjoint edges:
on the differing bits the endpoints are prescribed, and on the remaining
bits they agree. A closest pair therefore forces at least \(2^{k-d_*}\)
omissions. For \(A\) attaining \(d_*\) there are at most
\(\binom m2 2^{k-d_*}\) edges. Deleting one endpoint of each suffices.

A column of weight \(w\) separates \(w(m-w)\le\lfloor m^2/4\rfloor\)
pairs. Averaging over pairs gives the Plotkin bound \cite{P60},
\(d_*\le(1-\alpha_m)k\). For the
reverse estimate form a block consisting once of every column of weight
\(\lfloor m/2\rfloor\). By symmetry every pair is separated on exactly
the fraction \(1-\alpha_m\) of its columns. Repeating full blocks then
filling the bounded remainder gives \(d_*\ge(1-\alpha_m)k-O_m(1)\).
\end{proof}
For \(m=2\) antipodal words give deficit one from one test. For each fixed
\(m\ge3\), \(\alpha_m>0\), giving the exponential loss compared with the
polynomial loss for adaptive tests.

\clearpage
\section{Atlas of optimal strategies}\label{sec:atlas}
{\small Each row is one item. Column lengths are schematic. Letters name
widths locally: \(\mathrm{ABC}=A+B+C\), and two distinct blocks can both
have width \(D\). Take the listed rows and leftmost widths in the current
rectangles, top to bottom. Each successive test follows an answer of one.
Rectangles of the same color must be searched together. Grey marks a
single rectangle. Line styles distinguish tests. Brackets mark widths
in a tight bound for this pattern. The hatching selects the five stars of
the eight-item proof. Panel ranges give budgets with positive displayed
widths; Table~\ref{table:1} covers the earlier values.\par}

\begingroup
\setlength{\parskip}{0pt}
\setlength{\topsep}{3pt}
\setlength{\medskipamount}{3pt}
\setlength{\smallskipamount}{2pt}
\input{strategy-atlas}
\endgroup
\clearpage
\bibliographystyle{plain}
\bibliography{fixed-shore}

%% file: layout.tex
\usepackage[a4paper,margin=26mm,top=23mm,bottom=24mm,headheight=15pt]{geometry}
\usepackage{amsmath,amssymb,amsthm,booktabs,graphicx,caption,tikz}
\usepackage{lmodern,microtype,fancyhdr,needspace}
\usetikzlibrary{patterns,decorations.pathreplacing}
\usepackage[unicode,hidelinks]{hyperref}
\usepackage{bookmark}
\newtheorem{theorem}{Theorem}
\newtheorem{lemma}[theorem]{Lemma}

\newcommand{\papername}{Quantitative group testing with a fixed shore: sharp asymptotics and exact formulas}
\newcommand{\papersignature}{Fedor Karpelevitch\\[-1pt]{\small Independent researcher \(\cdot\) \href{mailto:fedor@karpelevitch.net}{fedor@karpelevitch.net}}}

%% file: width-changes.tex
\def\eightwidthchanges{1/0,2/+3,3/0,4/0,5/0,6/0,7/0,8/-3}

%% file: initial-values.tex
\begin{table}[htbp]
\centering\small
\setlength{\tabcolsep}{4pt}
\caption{Exact frontiers $n(k,m)$ for $1\le m\le8$ and $0\le k\le20$. Zero means that no nonempty opposite shore is searchable. The formulas in \S\ref{sec:table} give all continuations.}
\label{table:1}
\begin{tabular}{r*{8}{r}}
\toprule
$k\backslash m$ & 1 & 2 & 3 & 4 & 5 & 6 & 7 & 8 \\
\midrule
0 & $1$ & $0$ & $0$ & $0$ & $0$ & $0$ & $0$ & $0$ \\
1 & $2$ & $1$ & $0$ & $0$ & $0$ & $0$ & $0$ & $0$ \\
2 & $4$ & $3$ & $2$ & $1$ & $0$ & $0$ & $0$ & $0$ \\
3 & $8$ & $7$ & $5$ & $4$ & $3$ & $2$ & $2$ & $1$ \\
4 & $16$ & $15$ & $12$ & $10$ & $9$ & $7$ & $6$ & $5$ \\
5 & $32$ & $31$ & $27$ & $24$ & $22$ & $19$ & $17$ & $15$ \\
6 & $64$ & $63$ & $58$ & $54$ & $50$ & $46$ & $42$ & $38$ \\
7 & $128$ & $127$ & $121$ & $116$ & $109$ & $104$ & $97$ & $91$ \\
8 & $256$ & $255$ & $248$ & $242$ & $231$ & $225$ & $214$ & $206$ \\
9 & $512$ & $511$ & $503$ & $496$ & $481$ & $473$ & $457$ & $447$ \\
10 & $1\,024$ & $1\,023$ & $1\,014$ & $1\,006$ & $985$ & $973$ & $954$ & $939$ \\
11 & $2\,048$ & $2\,047$ & $2\,037$ & $2\,028$ & $2\,001$ & $1\,983$ & $1\,960$ & $1\,938$ \\
12 & $4\,096$ & $4\,095$ & $4\,084$ & $4\,074$ & $4\,039$ & $4\,015$ & $3\,988$ & $3\,956$ \\
13 & $8\,192$ & $8\,191$ & $8\,179$ & $8\,168$ & $8\,124$ & $8\,093$ & $8\,062$ & $8\,018$ \\
14 & $16\,384$ & $16\,383$ & $16\,370$ & $16\,358$ & $16\,304$ & $16\,265$ & $16\,228$ & $16\,172$ \\
15 & $32\,768$ & $32\,767$ & $32\,753$ & $32\,740$ & $32\,675$ & $32\,627$ & $32\,581$ & $32\,514$ \\
16 & $65\,536$ & $65\,535$ & $65\,520$ & $65\,506$ & $65\,429$ & $65\,371$ & $65\,315$ & $65\,237$ \\
17 & $131\,072$ & $131\,071$ & $131\,055$ & $131\,040$ & $130\,950$ & $130\,881$ & $130\,814$ & $130\,723$ \\
18 & $262\,144$ & $262\,143$ & $262\,126$ & $262\,110$ & $262\,006$ & $261\,925$ & $261\,846$ & $261\,741$ \\
19 & $524\,288$ & $524\,287$ & $524\,269$ & $524\,252$ & $524\,133$ & $524\,039$ & $523\,947$ & $523\,828$ \\
20 & $1\,048\,576$ & $1\,048\,575$ & $1\,048\,556$ & $1\,048\,538$ & $1\,048\,403$ & $1\,048\,295$ & $1\,048\,189$ & $1\,048\,054$ \\
\bottomrule
\end{tabular}
\end{table}